\documentclass[10pt]{article}
\usepackage{amsfonts}
\usepackage{amsmath}
\usepackage{amssymb}
\usepackage{color}

\def\eps{\varepsilon}

\def\R{\mathbb{R}}
\def\N{\mathbb{N}}

\def\Z{\mathbb{Z}}

\def\bo{\nl\phantom{a}\hfill $\Box$\nl}

\def\nl{\newline}

\begin{document}

\begin{center}
{\bf {\Large Classifying Functions via growth rates of repeated iterations\footnote{To appear in Fundamenta Mathematicae.}}} 
\vspace{0.2in}\\ Titus Hilberdink\\
Nanjing University of Information Science and Technology (Reading Academy)
219 Ningliu Road
Nanjing, China\\
E-mail: t.w.hilberdink@nuist.edu.cn, 
t.hilberdink@reading.ac.uk\vspace{0.15in}\\
\end{center}

\indent
\begin{abstract} In this paper we develop a classification of real functions based on growth rates of repeated iteration. We show how functions are naturally distinguishable when considering inverses of repeated iterations. For example, $n+2\to 2n\to 2^n\to 2^{\cdot^{\cdot^2}}$ ($n$-times) etc. and their inverse functions $x-2, x/2, \log x/\log 2,$ etc. Based on this idea and some regularity conditions we define classes of functions, with $x+2$, $2x$, $2^x$ in the first three classes. 

We prove various properties of these classes which reveal their nature, including a `uniqueness' property.  We exhibit examples of functions lying between consecutive classes and indicate
how this implies these gaps are very `large'. Indeed, we suspect the existence of a {\em continuum} of such classes.\nl

\vspace{0.3in}\end{abstract}

\noindent
{\bf Introduction}\nl
Consider a strictly increasing function $f_0:\N\to\N$ for which $f_0(n)>n+1$. Let $f_1(n) = f_0^n(1)$ (where $f^n$ denotes the $n^{\rm th}$-iterate of $f$) and more generally $f_{k+1}(n) = f_k^n(1)$ for each $k\in\N_0$. Then $f_1(n)>f_0(n)$ for $n>1$, and so $f_{k+1}(n)>f_k(n)$ for each $k$. Trivially then, each $f_k(n)\to\infty$ as $n\to\infty$, faster than $f_{k-1}(n)$. How fast? For example, with $f_0(n)=n+2$, we have
\[ f_1(n)=2n+1, \quad f_2(n)=2^{n+1}-1, \quad f_3(n) = 2^{\cdot^{\cdot^2}}-1\mbox{ ($n+1$-times)}\]
after which we run out of any standard notation. More interestingly, given a different starting function, say $g_0$, with $g_k$ defined in an analogous manner to $f_k$ above, how do $f_k$ and $g_k$ compare? In this article we attempt to answer this question. In doing so, we shall see that it leads to a natural classification of growth rates of real functions.

As a first step to distinguish between the behaviours at infinity, extend the functions to the reals and consider their inverses. Note that each $f_k$ extends to a strictly increasing continuous function, which has an inverse, say $F_k = f_k^{-1}$. It follows from $f_{k+1}(n+1)=f_k(f_{k+1}(n))$ that 
\[  F_{k+1}(f_k(x)) = F_{k+1}(x)+1\tag{1.1}\]
whenever $x=f_{k+1}(n)$. In fact we can always choose $F_{k+1}$ so that it satisfies (1.1) for {\em all} $x$ large enough; i.e. $F_{k+1}$ is an {\em Abel function} (see the Appendix, section 3) of $f_k$. We now drop the requirement that $f$ maps $\N$ to $\N$. 
\nl

\noindent
{\em The functions $\Xi_k$.}\, In order to investigate two rival `systems' (i.e. $f_k$ and $g_k$ as above), it is useful to have one such particular sequence of functions. We find it more convenient to do this for the inverses (i.e. the $F_k$) and we shall denote them by $\Xi_k$. We start with 
\[ \Xi_0(x)=x-e,\quad\Xi_1(x)=\frac{x}{e}, \quad\Xi_2(x)=\log x,\]
noting that (1.1) holds for $k=0,1$. For $\Xi_3$ we use the function $\Xi$ as defined in \cite{TH1} by $\Xi(x)=x$ on $[0,1)$ and extended to $\R$ via $\Xi(e^x)=\Xi(x)+1$. As such, $\Xi$ is a continuously differentiable and strictly increasing Abel function of exp. Note that $\Xi(x)=o(\log_k x)$ for every $k$, so tends to infinity very slowly. For $\Xi_k$ with $k\ge 4$ we do exactly the same -- take $\Xi_k$  to be a continuous strictly increasing Abel function of $\Xi_{k-1}^{-1}$. In fact, we shall require them to have some further regularity, namely that $\frac{\Xi_k(x)}{x}$ is decreasing. This can always be achieved, as we show in the appendix.
The functions we will be discussing will lie in\footnote{As in \cite{TH1}, $f\in{\rm SIC}_{\infty}$ means $f$ is strictly increasing, continuous, and tends to infinity.}
\[ \mathcal{F}_C = \{f\in {\rm SIC}_{\infty}\, :\, \Xi_k< f< \Xi_k^{-1}\quad\mbox{for some $k$}\}. \]
We call $\mathcal{F}_C$ the space of {\em functions of finite class}. It is independent of the choice of $\Xi_k$, in the sense that, starting with $f_0\in\mathcal{F}_C$, we find $\{f\in {\rm SIC}_{\infty}\, :\, F_k< f< F_k^{-1}\quad\mbox{for some $k$}\}=\mathcal{F}_C$.
\nl

Now, to study a sequence $f_k$ (or $F_k$) with a different starting point, it is more convenient to consider $F_{k+1}$ satisfying (1.1) with an extra $o(1)$ term, namely
\[ F_{k+1}(f_k(x)) = F_{k+1}(x)+1+o(1)\tag{1.2}\]
This notion was defined in \cite{TH1} and used to develop `orders of growth' of real functions (defined on a neighbourhood of infinity) with respect to another function.
To recall, a function $f$ has {\em order $\lambda\in\R$ with respect to a function} $F\in {\rm SIC}_{\infty}$ if
\[  F(f(x))-F(x)\to\lambda \quad\mbox{ as $x\to\infty$.} \]
This was denoted by $O_F(f)=\lambda$. Thus (1.2) says $O_{F_{k+1}}(f_k)=1$, or\footnote{We recall that $O_F(f\circ g) = O_F(f)+O_F(g)$ and $O_F(f^{-1})=-O_F(f)$ if $f^{-1}$ exists. So $O_F(f^k)=k$ for every integer $k$, where $f^k$ is the $k^{\rm th}$-iterate of $f$.} $O_{F_{k+1}}(F_k)=-1$.

For the $\Xi_n$, we started with $x+e$, but now we start with a general function $f$. We shall assume that 
$f\in{\rm SIC}_{\infty}$ is of finite class and $f(x)\geq x+1+\delta$ for some $\delta>0$. We can form a sequence of ${\rm SIC}_{\infty}$-functions $\{F_n\}_{n\geq 0}$ satisfying
\[ F_0 =f^{-1}\quad\mbox{ and }\quad  O_{F_{n+1}}(F_n)=-1,\tag{1.3} \]
for every $n\geq 0$. Of course the $F_n$ are not unique but, as we shall see later, given some extra regularity, they are unique up to asymptotic equivalence.
Let us consider a few examples and compare the resulting $F_k$. \nl

(i) \,  $f(x)=x+2$. Taking $F_1(x)=\frac{1}{2}x$, we find that (1.2) holds (without error term and $f_0=f$). Next we want $F_2$ such that $O_{F_2}(F_1)=-1$. Take $F_2(x)=\frac{\log x}{\log 2}$ so again (1.2) holds without error. For $F_3$, we can take $\Xi$, since $c\log x$ (with $c>0$) has order $-1$ with respect to $\Xi$. This is where the $o(1)$ comes into play.  
We see that after three steps we end up with $\Xi$, and we can take $F_n = \Xi_n$ for $n\geq 4$. \nl

(ii)\, $f(x)=x^2$. Now for $F_1$ we take $\frac{\log\log x}{\log 2}$ since then $F_1(f(x))-F_1(x) = 1$.
For the next step we can take $F_2(x)=\frac{1}{2}\Xi(x)$ since
\[ \Xi(F_1(x))-\Xi(x) = \Xi(\log\log x)-\Xi(x)+o(1) = -2+o(1),\]
giving $O_{F_2}(F_1)=-1$ as desired. Next we can take $F_3=\Xi_4$ (since $O_{\Xi_4}(F_2)=-1$), and so we can take $F_n=\Xi_{n+1}$ for $n\geq 4$. \nl

A few more examples are contained in the following table.\footnote{Here $a$ is a constant greater than 1. Also, by $f\simeq g+h$ we mean $f-g\sim h$.}
\[  \begin{array}{|c|cc@{\quad\vdots\quad}ccc@{\quad\vdots\quad}cc@{\quad\vdots\quad}c|}
\hline
f_0 & x+a  & x+\sqrt{x} & x+\frac{x}{\log x} & ax & x^a & e^x & e^{e^x} & \Xi^{-1}(x)\\
\hline
F_0 & x-a  & \simeq x-\sqrt{x} & \simeq x-\frac{x}{\log x} & \frac{x}{a} & x^{1/a} & \log x & \log\log x & \Xi(x)\\
F_1 & \frac{x}{a} & 2\sqrt{x} & \frac{1}{2}(\log x)^2 & \frac{\log x}{\log a} & \frac{\log\log x}{\log a} & \Xi(x) & \frac{1}{2}\Xi(x) & \Xi_4(x) \\
F_2 & \frac{\log x}{\log a} & \frac{\log\log x}{\log 2} & \Xi(x) & \Xi(x) & \frac{1}{2}\Xi(x) & \Xi_4(x) & \Xi_4(x) & \Xi_5(x)  \\
F_3 & \Xi(x) & \frac{1}{2}\Xi(x) & \Xi_4(x) & \Xi_4(x) & \Xi_4(x) & \Xi_5(x) & \Xi_5(x)& \Xi_6(x) \\
F_4 & \Xi_4(x) & \Xi_4(x) & \Xi_5(x) & \Xi_5(x) & \Xi_5(x) & \Xi_6(x)& \Xi_6(x) & \Xi_7(x) \\
\hline
\end{array} 
 \]
We notice that the functions split up into different behaviour. This observation leads naturally to the following definition.
\nl

\noindent
{\bf Definition 1.0}\, Let $f\in {\rm SIC}_{\infty}$ for which $f(x)\geq x+1+\delta$ (for some $\delta>0$), and let $k\in \N_0$. We say 
$f$ is of {\em class} $k$ if for functions $F_n\in {\rm SIC}_{\infty}$ satisfying (1.3) (for all $n\geq 0$), we can take
\[ F_n= \Xi_{n+k} \]
for some $n\geq 1$ (and hence all larger $n$). This is well-defined for, as we shall see later, if $f$ is of class $k$ and $G_n$ is any other sequence of functions satisfying (1.3) (i.e. with $G_0=f^{-1}$ and $O_{G_{n+1}}(G_n)=-1$) then we must have $G_n\sim\Xi_{n+k}$, provided the $\Xi_n$ satisfy some regularity condition (see Proposition 1.1). 
\nl

\noindent
Thus, for example, $x+2$ is of class 0, $2x$ and $x^2$ are of class 1, $e^x$ of class 2, $\Xi_k^{-1}$ of class $k$. \nl

The above definition, although natural from the preceding discussion, has some drawbacks. In particular, it is not true that functions of a certain class are necessarily smaller than those of a higher class. The above table would seem to suggest this, but it is false. For example, there are functions $f$ of class 1 for which $f(n)=n+1$ for infinitely many $n\in\N$ (see the Appendix).

In order to have the property that functions of a higher class are larger than those of a lower class, we need to add extra `regularity' conditions on our functions. In particular, the functions will be taken to be differentiable and their derivatives ``well behaved'', like the functions in our table. This leads to a more refined notion of classes, $C_n$ ($n\geq 0$), as defined in section 2. Roughly speaking, we construct a sequence of classes $C_n$ such that functions in the $n^{\rm th}$ class are `close' to $\Xi_n^{-1}$. More specifically, for functions $f$ in the $n^{\rm th}$-class, there exists an $F$ with $F^{-1}$ in the $n+1^{\rm th}$ class such that $O_F(f)=1$. 

There is an interesting but somewhat superficial resemblance to the classes of recursive functions defined by Grzegorczyk \cite{Grz} (see also \cite{Rit}) --- the so-called {\em Grzegorczyk hierarchy}. The difference is that these constructions are based on recursive complexity, whereas ours is based on rates of growth. The similarity is that the largest functions in the respective $k^{\rm th}$-classes have somewhat similar growth. 

This article is really an attempt to try to understand the nature of these classes; i.e. which functions lie in which class, which functions, say between $x^2$ and $e^x$ lie {\em between} classes 1 and 2, etc. As we shall see, all Hardy $L$-functions\footnote{Those functions obtained from a finite number of applications of $+,-$ and $\log$ and $\exp$ on the function $x$ and real constants (see \cite{H}).} greater than $x+a$ (with $a>1$) lie in $C_0,C_1,C_2$. It is not intuitively clear that any functions lie between consecutive classes. But such functions do exist (in contrast to the classes of recursive functions), though they are hard to find. Their growth rate is difficult to describe (other than saying they lie between, say, $C_1$ and $C_2$ functions) as they have no tangible connection to any known function. We discuss how their existence leads on to the surmise that there should exist a continuum of classes. This indicates there is a huge chasm between consecutive classes with only an infinitesimally small proportion of growth rates found within them. We cannot quite prove this but already show that there are big gaps between classes.

We also prove a `stability' or `uniqueness' result for these classes; namely that if they are defined using a different starting sequence, the resulting classes are either completely distinct or identical to the original ones.
\nl

There are two ways we can proceed, depending on whether we assume the existence of a `best behaved' set of functions $\{\Xi_n\}$, or not\footnote{The point is that there are many different functions satisfying (1.1) or (1.2) and it is not clear which particular choice is ``best'' in some sense.}. In the former case, we use a regularity condition based on $\sim$, while in the latter case is based on $\asymp$. We shall proceed as in the former case. In the Appendix, we show how the definitions and results adjust for the latter case.
\nl

To summarise briefly the rest of the paper; in \S1 below, we introduce some regularity on functions and prove various preliminary results which we shall need in order to define the classes of functions in \S2, and to be able to prove their basic properties. In particular, we prove in Theorem 2.3 and Corollary 2.4 that the classes are hierarchical (functions in a given class are larger than those in a previous class), a uniqueness result (Theorem 2.6), and the existence of functions between classes (Theorem 2.8).   

In \S3, we discuss further questions and open problems; in particular we link the idea of a continuum of classes to a problem regarding continuous extensions of the Ackermann function.

In the appendix, we add some background theory (eg on regular variation, Abel functional equation) and other results which give more context.  
\nl

\noindent
{\bf Notation.}\, Throughout this article, we shall use the same notation as that used in \cite{TH1}. 
Unless stated otherwise, all functions are considered to be defined on a neighbourhood 
of infinity. By $f<g$, we mean $f(x)<g(x)$ for all $x$ sufficiently large. The relations $f\sim g$ and $f=o(g)$ are defined as usual, namely: $f(x)/g(x)\to 1$ and 0 as $x\to\infty$, respectively. The relation $f\asymp g$ means there exist $a,b>0$ such that $a<f/g<b$ on a neighbourhood of infinity. Also $f\prec g$ means $f=o(g)$ while $f\succ g$ is just $g\prec f$.

We write $f^n$ for the $n^{\rm th}$-iterate, while for $\exp$ and $\log$ we use $e_n^{(\cdot)}$ and $\log_n {(\cdot)}$ respectively. 
Also
\begin{align*}
{\rm SIC}_{\infty} &= \{ f: f\mbox{ is continuous and strictly increasing and $f(x)\to\infty$
as $x\to\infty\}$,}\\
{\rm D}_{\infty}^+ &=\{f\in {\rm SIC}_{\infty}:f\mbox{ is continuously differentiable and 
$f^{\prime}>0\}$.}
\end{align*}
These are groups with respect to composition if we identify functions which agree on a neighbourhood of infinity.\nl

\noindent
{\bf {\large \S1. Regularity conditions and $\mathcal{B}_F$}}\nl
We shall require functions to satisfy various regularity conditions. 
For $F\in {\rm SIC}_{\infty}$, 
\[   F(x+o(x))=  F(x)+o(1)\tag{$R_0$} \]
For $F\in {\rm D}_{\infty}^+$,
\begin{align*}
  F^{\prime}(x+\lambda) &\sim F^{\prime}(x)\quad\mbox{ uniformly for $\lambda$ in bounded subsets of $\R$}  \tag{$R_1$} \\
  F^{\prime}(x+o(x)) &\sim F^{\prime}(x)\tag{$R_2$} \\
 F^{\prime}(\lambda x) &\sim \frac{1}{\lambda}F^{\prime}(x)\quad\mbox{ for every $\lambda>0$}.\tag{$R_3$} 
\end{align*}
Note that $(R_3)$ says $F^{\prime}$ is {\em regularly varying of index} $-1$ (see the appendix for a brief exposition on regular variation). 
Also $(R_3)\implies (R_2) \implies (R_1)$, while $(R_3)\not\!\!\!\implies (R_0)$ (take $F(x)=(\log x)^2$). Indeed $(R_0)\implies F(x)=O(\log x)$. Conversely, if $F^\prime(x)=O(\frac1x)$ then $F$ satisfies $(R_0)$.

For example, $\log x$ and $\Xi(x)$ satisfy $(R_0)$ and $(R_3)$.\nl

\noindent
{\bf Proposition 1.1}\nl
{\em Let $f,g\in {\rm SIC}_{\infty}$ such that $f\sim g$, and suppose $O_F(f)=O_G(g)=-1$ for some $F,G\in {\rm SIC}_{\infty}$, where $G$ (or $F$) satisfies ($R_0$). Then $F\sim G$.}\nl

\noindent
{\em Proof.}\, It suffices to show $O_G(f)=-1$, since then $O_G(f^{-1})=1 = O_F(f^{-1})$, which implies $F\sim G$ by Proposition 1.7 of \cite{TH1}. But $f\sim g$, so $G(f) = G(g+o(g)) =  G(g) +o(1) = G-1+o(1)$. Thus $O_G(f)=-1$.
\bo

In Proposition A1 of the Appendix we show that the $\Xi_n$ can be chosen to satisfy $(R_0)$ for $n\geq 2$. Hence Definition 1.0 can be replaced by: {\em $f$ is of class $k$ if, for all $F_n\in {\rm SIC}_{\infty}$ satisfying (1.2) with $F_0=f^{-1}$, we have $F_n\sim\Xi_{n+k}$ for some $n\geq 1$ (and hence all larger $n$)}.\nl

\noindent
{\bf The set} $\mathcal{B}_F$.\, Given $f,F\in {\rm D}^+_\infty$ such that $F\circ f\sim F$, it is of course not generally true that $(F\circ f)^\prime\sim F^\prime$. We call $\mathcal{B}_F$ the set of those $f$ for which this does hold.\nl

\noindent
{\bf Definition 1.1}\, For $F\in {\rm D}_{\infty}^+$, define $\mathcal{B}_F$ by
\[ \mathcal{B}_F = \{ f\in {\rm D}_{\infty}^+: (F\circ f)^{\prime}\sim F^{\prime}\}. \]
Equivalently, writing $L=1/F^{\prime}$, this says $f\in\mathcal{B}_F$ if and only if
\[ f^{\prime}\sim\frac{L(f)}{L}.\tag{1.4} \]
Condition (1.4) may be regarded as a regularity condition on $f$ with respect to $F$. Note that $\lambda x\in\mathcal{B}_F$ for all $\lambda>0$ is equivalent to $F$ satisfying $(R_3)$.
For example, taking $F=\log,\log\log,$ and $\Xi$ in turn, we have:
\begin{align*}
&  f\in\mathcal{B}_{\log}\mbox{ if }\quad f^{\prime}(x) \sim \frac{f(x)}{x} \\
&  f\in\mathcal{B}_{\log\log}\mbox{ if }\quad f^{\prime}(x) \sim \frac{f(x)\log f(x)}{x\log x} \\
&  f\in\mathcal{B}_{\Xi}\mbox{ if } \quad f^{\prime}(x) \sim \frac{\chi(f(x))}{\chi(x)}.
\end{align*}
(Here we wrote $\Xi^\prime=1/\chi$; $\chi=1$ on $[0,1]$ and $\chi(x)=x\chi(\log x)$.) Note that $\mathcal{B}_{\log}\subset \mathcal{B}_{\log\log}\subset \mathcal{B}_{\Xi}$.

Note further that every Hardy $L$-function $f$ such that $f(x)\to\infty$ as $x\to\infty$ lies in $\mathcal{B}_{\Xi}$. For there exists $n,k\in\Z$ such that $\log_n f(x)\sim \log_{n-k} x$. We may assume $n\ge 1,k$. Hence, on differentiating (and using the fact that for Hardy functions $f,g$ tending to infinity such that $f\sim g$, we have $f^\prime\sim g^\prime$),
\[ \frac{f^\prime(x)}{f(x)\log f(x)\cdots \log_{n-1}f(x)}\sim \frac{1}{x\log x\cdots \log_{n-k-1}x}.\]
Thus
\[  \frac{f^\prime(x)}{\chi(f(x))} \sim \frac{f(x)\log f(x)\cdots \log_{n-1}f(x)}{\chi(f(x))x\log x\cdots \log_{n-k-1}x}
=\frac{\chi(\log_{n-k}x)}{\chi(\log_n f(x))\chi(x)}\sim\frac1{\chi(x)}. \]
Indeed, the same proof works for functions from Boshernitzan's extension $E$ of Hardy's class of $L$-functions (see \cite{Bosh}, \cite{Bosh2}). This extension includes solutions of certain algebraic differential equations and is closed under differentiation and integration. Shackell \cite{Shack1} showed, using ideas of Rosenlicht \cite{Ros}, that the growth properties of such functions are similar to the Hardy functions.

We prove a number of basic properties of $\mathcal{B}_F$.\nl

\noindent
{\bf Proposition 1.2}
\begin{enumerate}
\item $\mathcal{B}_F$ {\em is a group under composition} ({\em indeed it is a subgroup of} ${\rm D}_{\infty}^+$).
\item {\em If $F$ satisfies ($R_2$) and $f,g\in\mathcal{B}_F$, then $f\sim g$ implies $f^{\prime}\sim g^{\prime}$.  In particular, if $F$ satisfies $(R_3)$, then $f(x)\sim \lambda x$ implies $f^{\prime}\to\lambda$.
For $\lambda=0$ this says $f(x)=o(x)$ implies $f^{\prime}\to 0$.} 
\item {\em Let $G\in {\rm D}_{\infty}^+$ satisfy $(R_2)$ and suppose $F\in\mathcal{B}_G$. Then
\begin{description}
\item[(i)] $f\in\mathcal{B}_F\implies f\in\mathcal{B}_G$;
\item[(ii)] $f\in\mathcal{B}_G$ and $F\circ f\sim F$ $\implies f\in\mathcal{B}_F$;
\item[(iii)] for any $h\in {\rm D}_{\infty}^+$, $\mathcal{B}_h\subset \mathcal{B}_{G\circ h}$.
\end{description}
}
\item {\em If $F$ satisfies ($R_3$), then $f\in\mathcal{B}_F$ implies $\alpha f\in\mathcal{B}_F$ for every $\alpha>0$.}
\item {\em Let $F,G\in {\rm D}_{\infty}^+$. If $F^{\prime}\sim G^{\prime}$, then $\mathcal{B}_F=\mathcal{B}_G$.}
\end{enumerate}

\noindent
{\em Proof.}\, (a)\, If $f,g\in \mathcal{B}_F$, then $f\circ g\in\mathcal{B}_F$, since
\[ (F\circ f\circ g)^{\prime} = (F\circ f)^{\prime}(g)g^{\prime}\sim F^{\prime}(g)g^{\prime}=(F\circ g)^{\prime}\sim F^{\prime}. \]
Also $f^{-1}$ exists and 
\[ (F\circ f^{-1})^{\prime} = \frac{F^{\prime}(f^{-1})}{f^{\prime}(f^{-1})}\sim F^{\prime}. \]
Thus $f^{-1}\in\mathcal{B}_F$. The identity, $f(x)=x$, is trivially in $\mathcal{B}_F$.\nl

(b)\, Writing $F^{\prime}=1/L$, we have $L(x+o(x))\sim L(x)$. Then $f\sim g$ implies 
\[  f^{\prime}\sim\frac{L(f)}{L}\sim \frac{L(g)}{L}\sim g^{\prime}.  \]
The second part follows immediately for $\lambda>0$. For $\lambda=0$, $L(f(x))=L(o(x))=o(L(x))$, so that $f^{\prime}\to 0$.
\nl

\noindent
(c)(i),(ii)\, Note that $f\in\mathcal{B}_F$ implies $F\circ f\sim F$. Let $F\circ f\sim F$. Then
\begin{align*}
f\in\mathcal{B}_G & \Longleftrightarrow G^{\prime}(f)f^{\prime}\sim G^{\prime}\\
& \Longleftrightarrow G^{\prime}(F(f))F^{\prime}(f)f^{\prime}\sim G^{\prime}(F)F^{\prime}\tag{since $G^{\prime}\sim G^{\prime}(F)F^{\prime}$}\\
& \Longleftrightarrow G^{\prime}(F)F^{\prime}(f)f^{\prime}\sim G^{\prime}(F)F^{\prime}\tag{since $F(f)\sim F$ and $G$ satisfies $(R_2)$}\\
& \Longleftrightarrow F^{\prime}(f)f^{\prime}\sim F^{\prime}\\
& \Longleftrightarrow f\in\mathcal{B}_F
\end{align*} 
(iii)\, Let $f\in\mathcal{B}_h$. Then 
\[ \bigl((G\circ h)(f)\bigr)^\prime = \bigl(G(h\circ f)\bigr)^\prime = G^\prime(h\circ f)(h\circ f)^\prime \sim G^\prime(h)h^\prime = (G\circ h)^\prime \]
so that $f\in\mathcal{B}_{G\circ h}$.

\noindent
(d)\, We have $(F(\alpha f))^{\prime} = \alpha F^{\prime}(\alpha f)f^{\prime}\sim F^{\prime}(f)f^{\prime}\sim F^{\prime}$.\nl

\noindent
(e) Suppose $F^\prime\sim G^\prime$. For any $f\in{\rm D}_{\infty}^+$,
\[ (G\circ f)^{\prime} = G^{\prime}(f)f^{\prime}\sim F^{\prime}(f)f^{\prime}=(F\circ f)^{\prime}.\]
Thus $f\in \mathcal{B}_G\Leftrightarrow f\in \mathcal{B}_F$.
\bo

\noindent
{\bf Remark}\, In particular, (b) shows that asymptotic equivalences in $\mathcal{B}_F$ may be differentiated. This has some strong consequences. Suppose $f,g$ are two Abel functions of some given function (so $f\sim g$ automatically). If $f,g\in\mathcal{B}_F$ for some $F$, then $f^\prime\sim g^\prime$. But this implies $f=g\, +$ constant; i.e. membership of $\mathcal{B}_F$ implies {\em unique} Abel solutions. 

Also, from (c)(i) and (ii), given $G$ satisfes $(R_2)$ and $F\in\mathcal{B}_G$, we find that $\mathcal{B}_F = \{f \in \mathcal{B}_G: F\circ f\sim F\}$. \nl

For the next result, we note that if $f^\prime \to\lambda\ge 0$, then
\[  f(x+\eps(x)) = f(x) +(\lambda +o(1))\eps(x).\tag{1.5}\]
for every $\eps(x)=o(x)$. Also, if $\lambda>0$, then $(f^{-1})^\prime \to 1/\lambda$. 
\nl

\noindent
{\bf Theorem 1.3}\nl 
{\em Let $F\in {\rm D}_{\infty}^+$ satisfy $(R_1)$ and let $g\in\mathcal{B}_F$. Then, for all $f$ and all $\lambda\in\R$}
\[  F(f) - F \sim \lambda F^{\prime}(g)\quad\Longleftrightarrow\quad O_g(f) =\lambda.\]
{\em For $\lambda=0$, the left-hand side is to be interpreted as $F(f)-F=o(F^{\prime}(g))$.}\nl

\noindent
{\em Proof.}\, Note that the assumption on $F$ implies $F(x+\mu)-F(x)\sim\mu F^{\prime}(x)$ for $\mu=O(1)$ since $F(x+\mu)-F(x)=\mu F^{\prime}(x+\theta_x)$ (for some $\theta_x$ in 
between $0$ and $\mu$), which in turn is asymptotic to $\mu F^{\prime}(x)$.\nl

$(\Longrightarrow)$\,  Let $h=F\circ g\circ F^{-1}$, so that $h^{\prime}\to 1$.
Then 
\begin{align*}
F(g(f)) &= h(F(f)) = h\Bigl(F + (\lambda+o(1))F^{\prime}(g)\Bigr)\\
& = h(F) + (\lambda+o(1))F^{\prime}(g)\tag{by (1.5) for $h$}\\
&= F(g) + (\lambda+o(1))F^{\prime}(g)\\
&= F(g + \lambda+o(1)).
\end{align*}
Hence $O_g(f)=\lambda$.\nl

$(\Longleftarrow)$\, Conversely, suppose $O_g(f) =\lambda$. Then 
\begin{align*}
F(f) & = h^{-1}( F( g(f))) = h^{-1}\Bigl(F(g+\lambda+o(1))\Bigr)\\
&= h^{-1}\Bigl(F(g) + (\lambda+o(1))F^{\prime}(g)\Bigr)\\
&= h^{-1}(F(g)) + (\lambda+o(1))F^{\prime}(g)\tag{by (1.5) for $h^{-1}$}\\
& = F + (\lambda+o(1))F^{\prime}(g),
\end{align*}
as required.
\bo

\noindent
{\bf {\large \S2. Finite classes}}\nl

\noindent
{\bf 2.1 Definitions}\nl
We are now ready to define the classes of functions alluded to in the introduction. We could base it on the functions $\Xi_n$ but actually there is no reason to restrict ourselves to this particular sequence. We could equally well start with a different initial sequence of functions. \nl

Let $\{f_n\}_{n\geq 0}$ be a sequence of ${\rm D}_{\infty}^+$-functions satisfying:
\begin{enumerate}
\item $f_0(x)>x+1+\eps$ for some $\eps>0$;
\item $F_n = f_n^{-1}$ satisfies $(R_2)$ for $n=0,1$ and $(R_3)$ for $n\geq 2$;
\item for $n\geq 0$, $O_{F_{n+1}}(f_n) = 1$ and $f_n\in\mathcal{B}_{F_{n+1}}$.
\end{enumerate}
We shall call such a sequence a {\em class generating sequence}. We restrict our attention to those functions which behave regularly w.r.t. this sequence. Let 
\[ \mathcal{B}  = \bigcup_{n=0}^{\infty}\mathcal{B}_{F_n}. \]

\noindent
{\bf Definition 2.1}\, For each $n\geq 0$, let $C_n^{(0)}=\{ f_n\}$ be the set containing just the one function $f_n$.
Now suppose that we have defined $C_n^{(k)}$ for all $n\geq 0$ and some $k\geq 0$. Then let $C_n^{(k+1)}$ be the set
\[ C_n^{(k+1)} = \{f\in\mathcal{B}:O_F(f) = 1\quad\mbox{for some $F$ with $F^{-1}\in C_{n+1}^{(k)}$}\}.\]
It follows immediately that $C_n^{(k)}\subset C_n^{(k+1)}$ for all $n,k\geq 0$. Now define $C_n$ by 
\[ C_n = \bigcup_{k=0}^{\infty} C_n^{(k)}. \]
{\bf Remarks 2.2.}
\begin{enumerate}
\item Note that (a) and (c) readily imply $F_1(x)\leq c_1x$ for some $0<c_1<1$ and $F_2(x)\leq c_2\log x$ for some $c_2>0$. We can, without loss of generality, assume a stronger version of (b), namely that for $n\geq 2$, $F_n^{\prime}$ is differentiable and $F_n^{\prime\prime}(x)\sim -F_n^{\prime}(x)/x$. We prove this in the appendix (Proposition A.2). Also we show there that such sequences can easily be generated (Proposition A.3). 
\item An immediate consequence is that $F^{-1}\in C_{n+1}$ implies $f=F^{-1}(F+1)\in C_n$. (For $O_F(f)=1$ and $f\in\mathcal{B}_F$, $F\in\mathcal{B}$ implies $f\in\mathcal{B}$ by Proposition 1.2(c).)
\item The classes $C_n$ depend of course on the functions $f_n$. With $F_n=\Xi_n$ for $n\geq 0$, we obtain the original classes but with a difference: here the functions must all satisfy a certain regularity w.r.t. the $f_n$; i.e. they must belong to $\mathcal{B}$. Thus in this case, $f\in C_n$ implies $f$ is of {\em class} $n$, but not conversely, since class $n$ functions need not even be differentiable. Furthermore, as we prove in Corollary 2.4 below, $C_n<C_{n+1}$ (i.e. all functions in $C_n$ are smaller than those in $C_{n+1}$ etc.) but there exist $f$ of class 0 such that $f\not < C_1$.
\end{enumerate}

Staying with $F_n=\Xi_n$, we see that, for any $a>1$, $a^x$ lies in $C_2^{(1)}$ (since $O_{\Xi}(a^x)=1$). Thus, using (b),  (with $F^{-1}(x)=a^x$) the function $f^{-1}(x)=F^{-1}(F(x)+1)=ax$ is in $C_1^{(2)}$. It follows that $x+a\in C_0^{(3)}$. In fact we can show much more, namely that every Hardy $L$-function greater than $x+a$ lies in $C_0$, $C_1$ or $C_2$ (see Proposition A4 in the Appendix). 
\nl 

\noindent
{\bf Theorem 2.3}\nl
{\em Let $f\in C_n$ and $g\in C_{n+1}$, where $n\geq 1$. Then}
\[ (g^{-1})^{\prime}=o((f^{-1})^{\prime}). \]
\noindent
{\em Proof.}\, By induction on $k$ in $C_n^{(k)}$. For $k=0$, the result is true for all $n\geq 1$, since $F_{n+1}(f_n(x))\sim F_{n+1}(x) = o(x)$ (as $n+1\geq 2 $). But 
$F_{n+1}\circ f_n \in\mathcal{B}$ so $F_{n+1}\circ f_n$ satisfies $(R_3)$ for
\[   \frac{(F_{n+1}\circ f_n)^\prime(\lambda x)}{(F_{n+1}\circ f_n)^\prime(x)}\sim \frac{F_{n+1}^\prime(\lambda x)}{F_{n+1}^\prime(x)}\sim\frac1{\lambda}. \]
By Proposition 1.2(b), $(F_{n+1}\circ f_n)^{\prime}\to 0$.

Suppose now the result is true for all $0\leq m\leq k$ for some $k$, and all $n\geq 1$. Let $f\in C_n^{(k+1)}$ and $g\in C_{n+1}^{(k+1)}$. Then there exist $F,G$ such that
$F^{-1}\in C_{n+1}^{(k)}$ and $G^{-1}\in C_{n+2}^{(k)}$ and 
\[ O_F(f)=O_G(g)=1. \]
By hypothesis, $G^{\prime}=o(F^{\prime})$; i.e. $G=h(F)$ where $h^{\prime}\to 0$. Hence, 
\[ G(f) = h(F(f)) = h(F+1+o(1))= h(F) +o(1)= G+o(1).\tag{2.1} \]
i.e. $O_G(f)=0$. In particular, since $f_n\in C_n^{(k+1)}$, we have $O_G(f_n)=0$. Also, from above $O_G(g^{-1}\circ f)=-1$. Now, since $n\geq 1$, we have $f_n(x)\geq (1+\eps)x$ for some $\eps>0$. Hence $O_G((1+\eps)x)=0$, and it follows that for every $\lambda>0$ $O_G(\lambda x)=0$. Thus $g^{-1}(f(x)) = o(x)$ as $x\to\infty$.
But $g^{-1}(f)\in\mathcal{B}$, so by Proposition 1.2(b), $(g^{-1}(f))^{\prime}\to 0$; i.e. $(g^{-1})^{\prime}=o((f^{-1})^{\prime})$. The inductive step is proved.
\bo

\noindent
{\bf Corollary 2.4}
\begin{enumerate}
\item {\em Let $f\in C_n$ for some $n\geq 0$. Then $O_{g^{-1}}(f)=0$ for all $g\in C_{n+2}$.}
\item {\em  For $n\geq 0$, $C_n<C_{n+1}$; i.e. if $f\in C_n$ and $g\in C_{n+1}$, then $f<g$.}
\item {\em Let  $n\geq 3$. Every $F$ with $F^{-1}\in C_n$ satisfies $(R_0)$ and $(R_3)$}. 
\item {\em Let $n\geq 2$. If $f\in C_n$ and $g\in\mathcal{B}$ such that $f\sim g$ or $f^{-1}\sim g^{-1}$, then $g\in C_n$.} 
\end{enumerate}
{\em Proof.}\, (a)\, $f\in C_n$, so there exists $F$ with $F^{-1}\in C_{n+1}$ such that $O_F(f)=1$. Let $G=g^{-1}$. By Theorem 2.3, it follows that $G=h\circ F$ where 
$h^{\prime}\to 0$. Hence $O_G(f)=0$ exactly as in (2.1).\nl

\noindent
(b)\, Let $f\in C_n$ and $g\in C_{n+1}$. Then there exists $G$ with $G^{-1}\in C_{n+2}$ such that $O_G(g)=1$. But $O_G(f)=0<O_G(g)$ by (a). Thus $f<g$.\nl

\noindent
(c)\, Since $f_1(x)\geq cx$ for some $c>1$ we have, by (a), $O_F(\lambda x)=0$ for all $\lambda>0$, whenever $F^{-1}\in C_n$ and $n\geq 3$. 
This implies $F$ satisfies $(R_0)$. Further, $F^{\prime}\sim \frac{F_k^{\prime}}{F^{\prime}_k(F)}$ for some $k$, so that 
\[ F^{\prime}(\lambda x)\sim \frac{F_k^{\prime}(\lambda x)}{F^{\prime}_k(F(\lambda x))}\sim \frac{F_k^{\prime}(x)}{\lambda F^{\prime}_k(F(x)+o(1))}\sim \frac{F_k^{\prime}(x)}{\lambda F^{\prime}_k(F(x))}\sim\frac{F^{\prime}(x)}{\lambda}.   \]
Hence $F$ satisfies $(R_3)$ also. (In fact here we just require $F(\lambda x)\sim F(x)$.)\nl

\noindent
(d) Let $F$ be such that $F^{-1}\in C_{n+1}$ and $O_F(f)=1$. Then $O_F(h)=0$ for all $h(x)\sim x$ from (c). Hence $O_F(g)=1$ (in either case) and $g\in C_n$. 
\bo

We also note that $f\in C_n^{(m)}$ implies $f(\lambda x)\in C_n^{(m)}$ for every $\lambda>0$ if $n\geq 2$ and $m\geq 1$.

%

Note that from Theorem 1.3, if $f\in C_n$ with $F$ such that $O_F(f)=1$ and $F^{-1}\in C_{n+1}$, then for all $k$ sufficiently large,
\[ F_k(f)-F_k\sim F_k^{\prime}(F). \]
In fact here we just need to take $k$ so that $F_k(F)\sim F_k$ (by Proposition 1.2(c)). This holds for $k=n+3$ at least by Corollary 2.4 (a) and, depending on $F$, may hold for $k=n+2$. Now $F_k^{\prime}$ is regularly varying of index $-1$. Thus $F_k^{\prime}\sim 1/H_k$ where $H_k\in {\rm D}_{\infty}^+$ is regularly varying of index 1 and moreover satisfies
\[ H_k^{\prime}(x)\sim \frac{H_k(x)}{x}.\tag{2.2} \]
(see section 2 of the Appendix). Combining with the above shows that $f$ and $F$ are related by $F_k(f)-F_k\sim 1/H_k(F)$. This suggests that, if $F^{-1}\in C_{n+1}$, then $F_k^{-1}(F_k + \frac{1}{H_k(F)})$ lies in $C_n$. \nl

\noindent
{\bf Theorem 2.5}\nl
{\em Let $n\geq 0$ and let $F^{-1}\in C_{n+1}$. Then for $k\geq n+3$ (and for $k=n+2$ if $F_{n+2}(F)\sim F_{n+2}$) 
\[ F_k^{-1}\Bigl(F_k + \frac{1}{H_k(F)}\Bigr)\in C_n,\tag{2.3} \]
where $H_k\in {\rm D}_{\infty}^+$ is any function asymptotic to $1/F_k^{\prime}$ satisfying (2.2)}.\nl

\noindent
{\em Proof.}\, Denote the function in (2.3) by $f$. Then $F_k(f)-F_k=\frac{1}{H_k(F)}\sim F_k^{\prime}(F)$, so that $O_F(f)=1$ by Theorem 1.3. It just remains to prove that $f\in \mathcal{B}$. It suffices to show $f\in\mathcal{B}_{F_k}$. But 
\[ (F_k\circ f)^{\prime} = F_k^{\prime}-\frac{H_k^{\prime}(F)F^{\prime}}{H_k(F)^2}. \]
The second term is, after (2.2), asymptotic to
\[ \frac{F^{\prime}}{F H_k(F)}\sim \frac{F^{\prime}F_k^{\prime}(F)}{F} = \frac{(F_k(F))^{\prime}}{F}\sim \frac{F_k^{\prime}}{F} = o(F_k^{\prime}). \]
Thus $(F_k\circ f)^{\prime}\sim F_k^{\prime}$, so $f\in\mathcal{B}_{F_k}$ and $f\in C_n$.\bo

\noindent
{\bf Examples}\nl
Now we give some more examples from the classes obtained from the original $\Xi_k$. With $n=0$ and $k=2$, Theorem 2.5 shows that if $F^{-1}\in C_1$ is such that $\log F(x) \sim \log x$ (i.e. $F(x)=x^{1+o(1)}$), then 
\[ \exp\Bigl\{ \log x + \frac{1}{F(x)}\Bigr\} = xe^{1/F(x)}\in C_0. \]
In the same way, but not directly following from Theorem 2.5; if $F^{-1}\in C_2$ and $\log F(x)> (\log x)^A$ for all $A$, then
\[ \exp\Bigl\{ \log x + \frac{1}{F(x)}\Bigr\} = xe^{1/F(x)}\in C_1. \]
(To show this one finds that, with $f$ denoting the above function, $O_G(f)=1$ where $G(x)=\frac{F(x)\chi(\log x)}{\chi(\log F(x))}$. By differentiating one finds $G^{\prime}(x)\sim \frac{F(x)}{x}$ which leads to $G\in\mathcal{B}$ and $G^{-1}\in C_2$.) Combining these gives (under the conditions)
\[ x+\frac{x}{F(x)}\in \left\{ \begin{array}{cl} C_0 & \mbox{ if }F^{-1}\in C_1\\ C_1 & \mbox{ if }F^{-1}\in C_2\end{array} \right. . \]
In this way, $x+x^\lambda\in C_0$ for $0<\lambda<1$, while $x+\frac{x}{(\log x)^\alpha}\in C_1$ for every $\alpha>0$. The function $\phi(x)=\Xi^{-1}(\Xi(x)+\frac12)$, which satisfies $\phi(\phi(x))=e^x$ lies in $C_2^{(2)}$ since, with $F=2\Xi$, $O_F(\phi)=1$ and $F^{-1}\in C_3^{(1)}$. The same is true for $\Xi^{-1}(\Xi+\beta)$ for any $\beta>0$, and even $\Xi^{-1}(\Xi+\frac1{\Xi})$.

We give some further examples to indicate functions in the $C_n$ can be quite `wobbly'. For instance, the function $f(x)=xh(x)$ where $h(x)= 3+\sin\Xi(x)$ lies in $C_1$, even though $f(x)\not\sim cx$ for any $c>0$, while of course $f(x)\asymp x$. Indeed, $f\in C_1^{(2)}$. For
\[ (\log f(x))^\prime = \frac{f^\prime(x)}{f(x)} = \frac{1}{x} + \frac{\cos\Xi(x)}{h(x)\chi(x)} \sim \frac{1}{x},\]
so $f\in\mathcal{B}$, and with $F(x) = \frac{\log x}{\log h(x)}$, we find that $O_F(f)=1$. A calculation shows $\frac{F^\prime(x)}{F(x)}\sim \frac{1}{x\log x}$ so $F\in\mathcal{B}$ and $O_{\Xi}(F)=-1$ since $F(x)\asymp \log x$. Thus $F^{-1}\in C_2^{(1)}$ and $f\in C_1^{(2)}$. 

In a similar way, the function $x^{h(x)}$ lies in $C_1^{(3)}$. This function wobbles between $x^2$ and $x^4$. 
\nl

\noindent
{\bf 2.2 Uniqueness of classes}\nl
Now we consider two separate systems of classes; suppose that $\{f_n\}$ and 
$\{g_n\}$ are two class generating sequences, both satisfying the same regularity conditions (i.e. for each $n$, $f_n\in \mathcal{B}_{G_k}$ for some $k$ where $G_k=g_k^{-1}$, and vice versa, $g_n\in\mathcal{B}_{F_k}$ for some $k$. This actually implies $\cup_{k\geq 0} \mathcal{B}_{F_k}=\cup_{k\geq 0} \mathcal{B}_{G_k}$.) We show below that if the resulting classes have some non-trivial intersection, then in fact they must be identical.
\nl

\noindent
{\bf Theorem 2.6}\nl
{\em Let $\{f_n\}$ and $\{g_n\}$ be two class generating sequences, generating classes $C_n$ and $D_n$ respectively. Suppose further that $\cup_{k\geq 0} \mathcal{B}_{F_k} = \cup_{k\geq 0} \mathcal{B}_{G_k}$ where $F_k=f_k^{-1}$ and $G_k=g_k^{-1}$. If $C_m\cap D_m\neq\emptyset$ for some $m$, then $C_n=D_n$ for all $n$.}\nl

\noindent
{\em Proof.}\, We shall show that $f_n\in D_n$ and $g_n\in C_n$ for all $n\geq 0$. For then, $C_n^{(0)}\subset D_n$ and $D_n^{(0)}\subset C_n$ and, by construction of $C_n$ and $D_n$, it follows that $C_n^{(k)}\subset D_n$ and $D_n^{(k)}\subset C_n$ for each $k\geq 0$ and $n\geq 0$. Hence $C_n\subset D_n$ and $D_n\subset C_n$, which implies $C_n=D_n$.

By assumption, there exists $k\geq 2$ such that $C_m^{(k)}\cap D_m^{(k)}\neq\emptyset$. Let $h\in C_m^{(k)}\cap D_m^{(k)}$. 
Then $\exists F,G$ with $F^{-1}\in C_{m+1}^{(k-1)}$ and $G^{-1}\in D_{m+1}^{(k-1)}$ such that
\[ O_F(h)=O_G(h)=1. \]
By Proposition 1.7 of \cite{TH1}, this implies $F(x)\sim G(x)$, and Corollary 2.4(d) then tells us that $F^{-1}\in D_{m+1}^{(k-1)}$ and $G^{-1}\in C_{m+1}^{(k-1)}$. Hence
$C_{m+1}^{(k-1)}\cap D_{m+1}^{(k-1)}\neq\emptyset$. Repeating this argument, we find that $C_{m+k-1}^{(1)}\cap D_{m+k-1}^{(1)}\neq\emptyset$. The same argument then 
shows that $f_{m+k},g_{m+k}\in C_{m+k}^{(1)}\cap D_{m+k}^{(1)}$. An easy induction shows $f_n,g_n\in C_n\cap D_n$ for $n\geq m+k$. (In fact, they lie in $C_n^{(1)}\cap 
D_n^{(1)}$, but we do not require this.) It then follows that $f_{n-1},g_{n-1}\in C_{n-1}\cap D_{n-1}$ for $n\geq m+k$ (since $O_{F_n}(f_{n-1})=1$ and $F_n^{-1}\in D_n$, and
similarly for $G_n$). Repeating this argument, we find that 
\[ f_n,g_n\in C_n\cap D_n\quad\mbox{ for $n\geq 0$.} \]
As shown above, this proves the theorem.
\bo

\noindent
{\bf 2.3 Functions with growth rates between classes}\nl
From the examples mentioned above, it may look like there are no functions between classes; i.e. that, for example, there is no function $f$ such that
$C_1<f<C_2$. This is false, as we shall presently show. 

We start with a more precise version of Theorem 2.5 which exhibits functions in $C_n^{(m+1)}$ which are greater (smaller) than everything in $C_n^{(m)}$.\nl

\noindent
{\bf Theorem 2.7}\nl
{\em Let $m,n\geq 1$. There exists $g,h\in C_n^{(m+1)}$ such that $g<C_n^{(m)}<h$; i.e. $g<f<h$ for every $f\in C_n^{(m)}$.

For $n=0$: let $m\ge 4$, then there exists $h\in C_0^{(m+1)}$ such that $h>C_0^{(m)}$; i.e. $h>f$ for every $f\in C_0^{(m)}$.}\nl

\noindent
{\em Proof.}\, We shall prove that such an $h$ exists in detail, then we show how the argument adapts easily to obtain $g$. 

For the proof, let $H_k$ be an asymptotic equivalent of $1/F^{\prime}_k$ which is in ${\rm D}_{\infty}^+$ and satisfies (2.2)\nl

\noindent
First consider $n\ge 1$. (We proceed by induction on $m$.) The case $m=1$. First note that for $f\in C_n^{(1)}$, $f\in\mathcal{B}$ and $O_{F_{n+1}}(f)=1$. This is equivalent to 
\[ F_{n+2}(f)-F_{n+2}\sim F_{n+2}^{\prime}(F_{n+1})\sim\frac{1}{H_{n+2}(F_{n+1})}. \]
Define $h$ by
\[ F_{n+2}(h)=F_{n+2} + \frac{2}{H_{n+2}(F_{n+1})}. \]
We claim that (i) $h>f$ for all $f\in C_n^{(1)}$, and (ii), $h \in C_n^{(2)}$.\nl

(i)\, From the above two formulas we have $F_{n+2}(f)-F_{n+2}<F_{n+2}(h)-F_{n+2}$ and $f<h$ follows.\nl

(ii)\, To show $h \in C_n^{(2)}$, we need to show $h\in\mathcal{B}$ and to find $F^{-1}\in C_{n+1}^{(1)}$ such that $O_F(h)=1$.
But as in the proof of Theorem 2.5, $h\in\mathcal{B}_{F_{n+2}}$ (the extra 2 is immaterial). Finally, with $F=\frac{1}{2}F_{n+1}$ we have $O_F(h)=1$ and $F^{-1}\in C_{n+1}^{(1)}$ since
\[ F_{n+2}(h)-F_{n+2}\sim\frac{1}{H_{n+2}(\frac{1}{2}F_{n+1})}  \]
(and using Theorem 1.3) and $(\frac{1}{2}F_{n+1})^{-1}\in C_{n+1}^{(1)}$ because $O_{F_{n+2}}(\frac{1}{2}x)=0$ (here we need $n\ge 1$). This proves the base case.

Now suppose the result holds for some $m\geq 1$ and all $n\geq 1$. Thus there exists $F$ such that $F^{-1}\in C_{n+1}^{(m+1)}$ and $F^{-1}>C_{n+1}^{(m)}$; i.e. for every $G$ with $G^{-1}\in C_{n+1}^{(m)}$, $F^{-1}>G^{-1}$. Equivalently, $F<G$ for all such $G$.

Define $h$ by 
\[  F_{n+3}(h)= F_{n+3} + \frac{2}{H_{n+3}(F)}. \]
Then, as in the $m=1$ case, it is easily seen that (i) $h>C_n^{(m+1)} $ and (ii) $h\in C_n^{(m+2)}$.  

For, given $f\in C_n^{(m+1)}$ with $G^{-1}\in C_{n+1}^{(m)}$ such that $O_G(f)=1$, we have 
\[  F_{n+3}(f)-F_{n+3}\sim \frac{1}{H_{n+3}(G)}  \]
(since $G\in \mathcal{B}_{F_{n+3}}$), while 
\[  F_{n+3}(h)-F_{n+3}= \frac{2}{H_{n+3}(F)}> \frac{2}{H_{n+3}(G)}>F_{n+3}(f)-F_{n+3}.  \]
Thus $h>f$. 

For (ii), we have $O_{\frac{1}{2}F}(h)=1$ and the inverse of $\frac{1}{2}F$ is in $C_{n+1}^{(m+1)}$.

\medskip

There is nothing special here about 2 in the definition of $h$; it can be replaced by any other constant greater than 1. To obtain a function $g\in C_n^{(m+1)}$ such that $g<C_n^{(m)}$, we just need to replace it by any constant between 0 and 1, say $\frac{1}{2}$. \nl

Now consider the case $n=0$. By the first part, we know there exists $F^{-1}\in C_1^{(m)}$ such that $F^{-1}>C_1^{(m-1)}$. Define $h$ by
\[  h = F_k^{-1}\Bigl(F_k+\frac{2}{H_k(F)}\Bigr)\]
(any $k\ge 3$. Then, as above, $h\in\mathcal{B}$, $h>C_0^{(m)}$ and $O_F(h)=2$. Hence $O_{\frac12F}(h)=2$. The result follows if we can show $(\frac12F)^{-1}\in C_1^{(m)}$; i.e $F^{-1}(2x)\in C_1^{(m)}$. Let $K^{-1}\in C_2^{(m-1)}$ be such that $O_K(F^{-1})=1$. It suffices to show $K(2x)=K(x)+o(1)$.

We require the following:\nl

\noindent
{\bf Lemma}\nl
{\em We have $(F_2\circ F_2)^{\prime}(x)=o(\frac1x)$.\nl

\noindent
Proof.}\, Since $(F_2\circ F_2)^{\prime}=F_2^\prime(F_2)F_2^{\prime}\sim\frac1{H_2 H_2(F_2)}$, so equivalently, we show $H_2(x) H_2(F_2)(x)\succ x$.

As $O_{F_2}(f_1)=1$ and $F_2$ satisfies $(R_3)$, we have\footnote{Here $F_2^2$ means $F_2\circ F_2$, not its square.} $O_{F_2^2}(f_1)=0$. But $f_1(x)>ax$ for some $a>1$. Thus $O_{F_2^2}(\lambda x)=0$ for all $\lambda>0$. Thus
\[  \int_{\frac{x}{2}}^x -\Bigl(\frac{1}{F_2^2}\Bigr)^\prime = \frac{1}{F_2^2(\frac{x}{2})} - \frac{1}{F_2^2(x)} = \frac{F_2^2(x)-F_2^2(\frac{x}{2})}{F_2^2(\frac{x}{2})F_2^2(x)}=o\Bigl(\frac{1}{F_2^2(x)^2}\Bigr).\]
But the integral is
\[  \int_{\frac{x}{2}}^x \frac{(F_2^2)^\prime}{(F_2^2)^2}\sim \int_{\frac{x}{2}}^x \frac1{H_2 H_2(F_2)(F_2^2)^2}\ge \frac{x}{2H_2(x) H_2(F_2(x))F_2^2(x)^2}.\]
Thus $H_2(x) H_2(F_2(x))\succ x$ follows.
\bo

\noindent
{\em Proof of Theorem 2.7 continued}\, Note that $f_2^2\in C_2^{(2)}$, for $O_{\frac12F_3}(f_2^2)=1$ and $F_3^{-1}(2x)\in C_3^{(1)}$. It follows (after Remark 2.2(b)) that
\[  \ell := f_2^2(F_2^2 + 1)\in C_1^{(3)}\]
and $F^{-1}>\ell$. But $K(F^{-n})\sim n\sim F_2^2(\ell^n)<F_2^2(F^{-n})$. Hence $K\ll F_2^2$. Thus
\[ K^\prime(x)\sim\frac{H_k(K(x))}{H_k(x)}\ll \frac{H_k(F_2^2(x))}{H_k(x)}\sim \frac{1}{H_k(x)H_k(F_2(x))}=o\Bigl(\frac1x\Bigr),\]
and, as a result, $K(2x)-K(x) = \int_x^{2x}K^\prime = o(1)$.
\bo

\noindent
{\bf Theorem 2.8}\nl
{\em Let $n\geq 0$. There exists $f$ such that $C_n<f<C_{n+1}$.}\nl

\noindent
{\em Proof.}\, For $n\ge 1$, we construct sequences of  functions $\{ g_m\}$ and $\{ h_m\}$ for $m\geq 1$ as follows. (For $n=0$ we start at $m=4$.) 

(i) \, Choose $g_1\in C_n^{(1)}$ and $h_1\in C_{n+1}^{(1)}$ such that $f_n<g_1< h_1<f_{n+1}$. Without loss of generality we may assume $g_1< h_1$ on $(0,\infty)$, for suppose it is true on $[x_0,\infty)$, then we may define $g_1(x)=\frac{x}{x_0}g_1(x_0)$ and $h_1(x)=\frac{x}{x_0}h_1(x_0)$ for $0< x< x_0$. 

(ii)\, Now suppose we have defined $g_k\in C_n^{(k)}$ and $h_k\in C_{n+1}^{(k)}$ for $k=1,\ldots, m$ so that 
\[  g_1< g_2< \ldots < g_m< h_m< \ldots < h_2< h_1 \qquad\mbox{ on $(0,\infty)$} \tag{2.4} \]
and $g_k> C_n^{(k-1)}$ and $h_k< C_{n+1}^{(k-1)}$ for each such $k$. 

By Theorem 2.7, there exist $g\in C_n^{(m+1)}$ and $h\in C_{n+1}^{(m+1)}$ such that $g> C_n^{(m)}$ and $h< C_{n+1}^{(m)}$.
In particular, $g_m<g<h<h_m$ on some interval $[x_1,\infty)$. Define $g_{m+1}=g$ and $h_{m+1}=h$ on $[x_1,\infty)$ and extend to $(0,x_1]$ so that $g_{m+1}<h_{m+1}$ on $(0,\infty)$, as we did for $g_1$ and $h_1$. Thus (2.4) holds for $m+1$. 

Since $g_m$ is monotonic and bounded, $\lim_{m\to\infty} g_m$ exists on $(0,\infty)$, and equals $f$ say.  
But then, for each $m$,
\[ C_n^{(m)}<g_{m+1}<f<h_{m+1}<C_{n+1}^{(m)} . \]
i.e. $C_n<f<C_{n+1}$.
\bo

\noindent
{\bf Heuristic discussion --- A Continuum of classes?}\nl
A function lying between consecutive classes, say, between $C_1$ and $C_2$, has a rather alien growth rate. Other than having a number of upper and lower bounds (due to the fact it lies between $C_1$ and $C_2$) there is nothing tangible to measure its growth rate by. This is contrast for example with the function $\phi(x)=\Xi^{-1}(\Xi(x)+\frac12)$ which satisfies $\phi(\phi(x))=e^x$. This function, though lying outside Hardy's class of $L$-functions (or indeed Boshernitzan's class $E$), lies in $C_2$ and can be ``measured'' via $\Xi$; namely, $\Xi(\phi)=\Xi+\frac12$. 

Let us now consider the likely consequences of Theorem 2.8. We have found a function $h$ such that $C_0<h<C_1$. From this, we can find a class generating sequence $\{h_n\}$ (with $h_0=h$) which generates classes $D_n$, say. If we can further ensure that $h_n\in\mathcal{B}$ for each $n$, then we should have $C_n<h_n<C_{n+1}$ for each $n$, which in turn results in $C_n<D_n<C_{n+1}$.  

To aid our intuition, let us write $D_n = C_{n+\frac{1}{2}}$. In the same way that a function can be found between two neighbouring classes, we can presumably do the same here and find a function -- and hence a class generating sequence, and hence classes -- between the $C_n$, $C_{n+\frac{1}{2}}$ and $C_{n+1}$. Call them $C_{n+\frac{1}{4}}$ and 
$C_{n+\frac{3}{4}}$ respectively. Repeating this process indicates that there are classes $C_{n+q}$ for any rational $q$ of the form 
$r/2^k$ with $r=0,1,\ldots, 2^k-1$ such that 
\[ C_{n+q}<C_{n+q^{\prime}} \]
for all such rationals $q,q^{\prime}$ with $q<q^{\prime}$. By taking limits of sequences of functions in $C_q$ with $q$ tending to any $\alpha\in [0,1]$, we should be able to find functions, and hence classes of functions, in between each of the $C_{n+q}$s. This suggests there exists a {\em continuum} of classes $C_{\alpha}$ 
$(\alpha\geq 0)$, such that\nl
 (i) $C_{\alpha}<C_{\beta}$ whenever $\alpha<\beta$ and\nl
 (ii) $f\in C_{\alpha}$ implies the existence of $F$ with $F^{-1}\in C_{\alpha+1}$ such that $O_F(f)=1$.

This indicates that the $C_n$ $(n\in\N_0)$ represent only an {\em infinitesimally small proportion} of the functions of finite class, and that there is a vast chasm of growth rates between consecutive classes.\nl

It would be most interesting if this can all be done in such a way that the (hypothesized) $C_{\alpha}$ are `complete' in the sense that there are no functions in between them; i.e.
\[ \not\!\!\exists\, f\quad\mbox{\em such that }\quad C_{\alpha}<f<C_{\alpha+\eps}\quad\mbox{\em for all $\eps>0$}.\tag{2.5} \] 

As suggested above, we expect a `vast chasm' between classes. We are not ready to prove the above continuum of classes but we can give some weaker results which already indicate there is a large gap between consecutive classes. \nl

\noindent
{\bf Proposition 2.9}\nl
{\em Let $k\ge 1$ and $f\in {\rm SIC}_\infty$ be such that $C_k<f<C_{k+1}$. Then, for any $F\in {\rm SIC}_\infty$ satisfying $O_F(f)=1$, we have $C_{k+1}<F^{-1}<C_{k+2}$.\nl

\noindent
Proof.}\, Let $G^{-1}\in C_{k+1}$ and $H^{-1}\in C_{k+2}$. Put $G_1 = \frac{1}{2}G$ and $H_1=2H$. Then $g:=G_1^{-1}(G_1+1)\in C_k$ and $h:=H_1^{-1}(H_1+1) \in C_{k+1}$. By assumption, $g<f<h$. Thus also $g^n(x)<f^n(x)<h^n(x)$ for all $n\in\N$ and $x\ge x_0$, sufficiently large. Thus 
\[ G(f^n(x))>G(g^n(x))\sim 2n \sim 2F(f^n(x)) \]
as $n\to\infty$, uniformly for $x$ in any bounded in interval $[x_0,A]$. Hence $G>F$. 

Similarly, 
\[  2H(f^n(x))<2H(h^n(x))\sim n \sim F(f^n(x)) \]
as $n\to\infty$, and $H<F$ follows. Thus $G^{-1}<F^{-1}<H^{-1}$, as required.
\bo

\noindent
{\bf Proposition 2.10}\nl
{\em Let $k\ge 1$ and $F^{-1}\in {\rm SIC}_\infty$ be such that $C_{k+1}<F^{-1}<C_{k+2}$. Put 
\[  f=F_m^{-1}\Bigl(F_m+\frac{1}{H_m(F)}\Bigr), \tag{2.6}\]
where $m$ sufficiently large. Then $C_k<f<C_{k+1}$. \nl

\noindent
Proof.}\, Let $g\in C_k$ and $h\in C_{k+1}$. Then there exist $G^{-1}\in C_{k+1}$ and $H^{-1}\in C_{k+2}$ such that $O_G(g) = O_H(h)=1$. As such, with $m\ge k+4$, 
\[ F_m(g) = F_m + \frac{1+o(1)}{H_m(G)}\quad\mbox{ and }\quad F_m(h) = F_m + \frac{1+o(1)}{H_m(H)}.  \]
Let $G_1 = \frac{1}{2}G$ and $H_1=2H$. Then $G_1^{-1}\in C_{k+1}$ and $H_1^{-1}\in C_{k+2}$ (here we need $k\ge 1$).
Thus $G^{-1}<F^{-1}<H^{-1}$ and so $H_m(G_1)>H_m(F)>H_m(H_1)$. But $H_m(G)\sim 2H_m(G_1)$ and $H_m(H)\sim\frac{1}{2}H_m(H_1)$. Thus
$F_m(g)-F_m<F_m(f)-F_m<F_m(h)-F_m$; i.e. $g<f<h$, as required.
\bo

\noindent
{\bf Remarks} \, (i)\, With a minor adjustment, the method in Proposition 2.10 shows that for any $c>0$, the function $F_m^{-1}(F_m+\frac{c}{H_m(F)})$ also lies between $C_k$ and $C_{k+1}$. (Take $G_1=\lambda G$ and $H_1=\mu H$ with $0<\lambda<c<\mu$.) Thus every function of the form
\[  F_m^{-1}\Bigl(F_m+\frac{\theta}{H_m(F)}\Bigr) \]
with $\theta\asymp 1$ and $m$ large enough lies between $C_k$ and $C_{k+1}$.

\noindent
(ii) It is not clear if the function defined by (2.6) is necessarily strictly increasing. If we want this, we could instead consider inverses: namely, define instead $f$ via 
\[ f^{-1}=F_m^{-1}\Bigl(F_m-\frac{1}{H_m(F)}\Bigr), \]
This function and its inverse is strictly increasing and $C_k<f<C_{k+1}$.
\nl

\noindent
{\bf A vast chasm between classes}\nl
 From Proposition 2.9, we can form a sequence  $\{f_n\}_{n\ge 1}$ such that $O_{F_{n+1}}(f_n) =1$ where $F_n=f_n^{-1}$, satisfying $C_n<f_n<C_{n+1}$. If we could show that the $f_n$ all lie in $\mathcal{B}$ and satisfy some smoothness conditions, then we could readily generate classes between $C_n$ and $C_{n+1}$. However, we do not know for certain this can be done (at this stage) but using Proposition 2.10, we can do something similar. Define $D_n^{(k)}$ recursively as follows: first let
\[ D_n^{(0)} = \{ f: f^{-1}=F_m^{-1}\Bigl(F_m-\frac{1+\delta}{H_m(F_{n+1})}\Bigr)\Bigr\},  \]
where $\delta$ is any continuous function decreasing to 0. By Proposition 2.10 (and Remark (iii) above) each $C_n< D_n^{(0)} < C_{n+1}$. Now, suppose $D_n^{(k)}$ has been defined for some $k\ge 0$ and all $n\ge 1$, let
\[ D_n^{(k+1)} = \Bigl\{ f: f^{-1} = F_m^{-1}\Bigl(F_m - \frac{1+\delta}{H_m(F)}\Bigr) \mbox{ where $\delta\searrow 0$ and } F^{-1}\in D_{n+1}^{(k)} \Bigl\}. \]
By induction we find that $C_n< D_n^{(k)} < C_{n+1}$ for all $k\ge 0$. Let $D_n = \cup_{k\ge 0} D_n^{(k)}$. Then $C_n< D_n < C_{n+1}$.

We could continue as in Theorem 2.7 showing the existence of functions $g_m\in C_n^{(m)}$, $h_m\in D_n^{(m)}$ such that $g_m>C_n^{(m-1)}$ and $h_m<D_n^{(m-1)}$ and
$g_m<g_{m+1}<h_{m+1}<h_m$ to obtain functions between the $C_n$ and $D_n$. However, we omit this for now in the hope that we obtain $f_n\in\mathcal{B}$ first, or for a different approach as discussed in the next section.
\nl

\noindent
{\bf {\large \S 3. Further thoughts and open questions}}\nl

\noindent
{\bf 3.1 An equivalence relation}\nl
An alternative way of looking at the problem highlighted in the heuristic discussion
is to go back to the original definition of classes (Definition 1.0). Given $f$ such that $f(x)>x+1+\delta$ for some $\delta>0$, we can always find a sequence of functions $\{f_n\}_{n\geq 0}$ with $f_0=f$ and $f_n\in {\rm SIC}_{\infty}$ satisfying $O_{F_{n+1}}(f_n)=1$ for each $n\geq 0$, where $F_n=f_n^{-1}$. We shall also insist that $F_n$ satisfies $(R_0)$ for $n\geq 2$. After Proposition 1.1 we note that for any two such sequences (with the same starting point $f$), say $\{f_n\}$ and $\{ \phi_n\}$ with $\Phi_n=\phi_n^{-1}$, we always have
\[ F_n\sim \Phi_n, \]
for all $n$. \nl

\noindent
{\bf Definition}\, For $f, g\in {\rm SIC}_{\infty}$ both larger than $x+1+\delta$ (some $\delta>0$), we say $f$ is {\em class-equivalent} to $g$, which we denote by 
\[ f\sim_c g, \]
if for sequences $\{ f_n\}$ and $\{ g_n\}$ as described above (with inverses $F_n$ and $G_n$ respectively), we have 
\[ F_n\sim G_n \]
for all $n$ sufficiently large.\nl

Thus $f\sim_c f$ always, and $\sim_c$ is clearly seen to be an equivalence relation. The connection with the previously defined classes arises from the fact that if $f,g\in C_k$, then $f\sim_c g$. This is true because the corresponding $F_n$ and $G_n$ eventually coincide with the underlying sequence defining the classes. Thus
\[ f,g\in C_k\Longrightarrow f\sim _c g. \]
Since $\sim_c$ is an equivalence relation, we have equivalence classes of functions. The question is now: amongst the functions of finite class, how many equivalence classes are there? Of course, the cardinality is at most that of the continuum as we are considering real continuous functions. The question is: {\em are there uncountably many, and is the cardinality that of the continuum?}\nl

\noindent
{\bf 3.2 The Ackermann function and generalizations}\nl
One way of (possibly) getting a handle on the hypothesized continuum of classes is by considering two-variable extensions of the {\em Ackermann} function (see \cite{Ack}).\nl 

There are several variants of the Ackermann function, all defined recursively and satisfying (3.1) below. We shall define it to be $A:\N^2_0\to \N$ as follows: 
$A(m,0)=2$, $A(0,n)=n+2$ and for $m,n\geq 0$,
\[ A(m+1, n+1) = A(m,A(m+1,n)).\tag{3.1} \]
It is easy to see that $A(1,n)=2n+2$ and $A(2,n)=2^{n+2}-2$.
\nl

\noindent
{\bf Generalizations}\nl
It is of interest to study extensions of the Ackermann function to real variables. For each fixed $m$, we can define an extension of $A(m,n)$ to $A(m,x)$ where $x$ ranges over real values in $[0,\infty)$, in such a way that each such function is strictly increasing and continuous and $G(m,x) := A^{-1}(m,x)$ satifies the Abel functional equation
\[ G(m+1, x) = G(m+1, G(m,x)) +1. \]
For example, we can define $G(0,x)=x-2$, $G(1,x)=\frac{1}{2}x-1$, and 
\[ G(2,x) = \frac{\log (x+2)}{\log 2} - 2. \]
In this way, the sequence of functions $\{A(m, \cdot)\}_{m\geq 0}$ can be seen to generate classes $C_m$. A suitable generalization in the second variable $m$ to real values, should give us an idea and possibly a definition for the continuum of classes we are seeking. So suppose we have extended $A(m,\cdot)$ to $A(\alpha,\cdot)$ (where $\alpha\in\R$ and $\alpha\geq 0$). Then for each $\alpha\in [0,1)$, $\{A(m+\alpha, \cdot)\}_{m\geq 0}$ generates classes, which we may denote by $C_{m+\alpha}$.  Of course, we want to make sure that each of the classes is separate and in particular that $C_{\alpha}<C_{\beta}$ whenever $\alpha<\beta$. For this reason we must ensure that any such generalization satisfies the condition
\[ G({\beta}, x) = o(G(\alpha, x)), \]
whenever $1\leq\alpha<\beta$. The question is whether this is possible, and moreover, whether it can be done in such a way that (2.5) holds for these classes.\nl

\noindent
{\bf 3.3. The operator $\mathcal{L}$}\nl
Finally we may view the problem in terms of an operator. Let ${\rm SIC}_\infty^+ = \{ f\in {\rm SIC}_\infty: f(x)>x\}$ and ${\rm SIC}_\infty^- = \{ f\in {\rm SIC}_\infty: f(x)<x\}$. These are closed under composition. Regarding inverses, note that $f\in {\rm SIC}_\infty^+\Leftrightarrow f^{-1}\in {\rm SIC}_\infty^-$. \nl

\noindent
{\bf Definition}\, Let $\mathcal{L}:{\rm SIC}_\infty^+\to{\rm SIC}_\infty^+$ denote the mapping $f\mapsto f(f^{-1}+1)$. \nl

Observe that $f^{-1}$ is an {\em Abel function} of $\mathcal{L}f$; i.e. $f^{-1}(\mathcal{L}f) = f^{-1}+1$. Equivalently, we can regard $\mathcal{L}^{-1}f$ as the set of inverse Abel functions of $f$. Observe that if $F^\prime(x+1)\sim F^\prime(x)$, then $f\in\mathcal{B}_F\implies \mathcal{L}f\in\mathcal{B}_F$. Hence $\mathcal{L}(C_{n+1})\subset C_n$ for each $n\ge 0$. 

Connecting to the discussion on the generalized Ackermann function, if we write $A_\alpha(x) = A(\alpha,x)$, then $\mathcal{L}A_{\alpha+1} = A_\alpha$. Thus, formally, we can view $A_\alpha$ as $A_\alpha = \mathcal{L}^{-\alpha}A_0$. In other words, we have to make sense of fractional iterates of this operator. The simplest non-trivial case is making sense of $\mathcal{L}^{\frac{1}{2}}$; i.e.{\em can we find an operator $\mathcal{H}$ on ${\rm SIC}_\infty^+$ such that $\mathcal{H}^2 = \mathcal{L}$}?\nl

\begin{center}
{\bf {\large Appendix}}
\end{center}

\noindent
{\bf 1. Examples of functions of class 0 which are not less than all class 1 functions}\nl
We give examples to show the existence of functions of a given class which are not smaller (respectively larger) than those of a higher (resp. lower) class.
\begin{enumerate}
\item Let $F$ be a ${\rm SIC}_{\infty}$ function for which $F(2^k)=2k$ and $F(2^k-1)=2k-1$ for each positive integer $k$. Note that for $2^k\leq x<2^{k+1}$, $2k\leq F(x)<2k+2$, so 
$F(x) = 2\frac{\log x}{\log 2}+O(1)$. 

Let $f=F^{-1}(F+1)$. Thus $O_F(f)=1$ and $O_{\Xi}(F)=-1$. Hence we can take $F_1=F$, $F_2=\Xi$, and $F_n=\Xi_{n+1}$ for $n\geq 2$. Thus $f$ is of class 1. 

But $F(2^k)=F(2^k-1)+1=F(f(2^k-1))$ so that $f(2^k-1)=2^k$, and similarly $f(2^k) = 2^{k+1}-1$; i.e. $f(x)=x+1$ and $f(y)=2y-1$ for some arbitrarily large $x$ and $y$.\nl

Note that we can easily adjust the example to make $f$ be as we large as we please while also keeping $f(x)=x+1$ for some arbitrarily large $x$. Just let $a_n$ be such that $a_{n+1}>a_n+1$ for all $n\ge 1$ and let $F\in {\rm SIC}_\infty$ satisfy $F(a_n)=2n$ and $F(a_n-1)=2n-1$. Putting $A(x)=\sum_{a_n\leq x} 1$, we see that $F(x) = 2A(x)+O(1)$. Thus $F^{-1}$ can be made to be of any given class, say class $k$ (eg. we can take $a_n=\Xi^{-1}_k(n)$). Then $f$ is of class $k-1$.

But $f(a_k-1)=a_k$, so $f(x)=x+1$ for some arbitrarily large $x$.
\item Let $a_k=2^{2^k}$. Let $F$ be a ${\rm SIC}_{\infty}$ function for which $F(2a_k)=F(a_k)+1$ and $F(n+2)=F(n)+1$ for each positive integer $n$ with $2a_k\leq n\leq a_{k+1}-2$.
Again, letting $f=F^{-1}(F+1)$, we find that $f(a_k)=2a_k$ while $f(n)=n+2$ for $2a_k\leq n\leq a_{k+1}-2$. Thus  $f(x)=x+2$ and $f(y)=2y$ for some arbitrarily large $x$ and $y$.

It remains to show that $f$ is of class 0. We do this by showing there exist $0<a<b<1$ such that $ax<F(x)<bx$. 

For each integer $r\geq 0$ such that $2a_k+2r\leq a_{k+1}$, we have $F(2a_k+2r)=F(2a_k)+r = F(a_k)+r+1$. In particular, $F(a_{k+1})-F(a_k)=\frac{1}{2}a_{k+1}-a_k+1$.
It quickly follows that $F(a_k)\sim\frac{1}{2}a_k$. Hence for all $k$ sufficiently large,
\[ F(2a_k+2r)=F(a_k)+r+1\leq a_k+r \]
for $0\leq 2r\leq a_{k+1}-2a_k$; i.e. $F(2n)\leq n$ for $2a_k\leq 2n\leq a_{k+1}$. But for $a_k\leq n<2a_k$, $F(n)\geq F(a_k)\sim\frac{1}{2}a_k\geq\frac{1}{4}n $. Since $F$ is increasing the result follows.

But now, any function $F_2\in {\rm SIC}_{\infty}$ for which $O_{F_2}(F)=-1$ necessarily has $F_2(x)\asymp \log x$. This is because $cx\leq F^{-1}(x)\leq dx$ for some $d>c>1$, so that $c^nx_0\leq F^{-n}(x_0)\leq d^nx_0$ for each fixed $x_0$ sufficiently large, and $F_2(F^{-n}(x_0))\sim n$. 

We may therefore take $F_3=\Xi$, $F_4=\Xi_4$ etc. Hence $f$ is of class 0.
\end{enumerate}

\noindent
{\bf 2. Regular Variation}\nl
A function $f:[A,\infty)\to\R$, defined on a neighbourhood of infinity, is {\em regularly varying with index} $\rho$ 
if it is locally bounded, measurable, and satisfies
\[ \lim_{x\to\infty} \frac{f(\lambda x)}{f(x)} = \lambda^{\rho} \quad\mbox{ for all $\lambda>0$}. \]
We note that this limit actually holds uniform in $\lambda$ on compact subsets of $(0,\infty)$. For the case $\rho=0$, $f$ is said to be {\em slowly varying}.

As such $f$ can be represented by
\[ f(x) = c(x)x^{\rho}\exp\Bigl\{ \int_A^x \frac{\eps(t)}{t}\, dt\Bigr\}, \]
where $c(x)\to c>0$ and $\eps(x)\to 0$ as $x\to\infty$. It follows that $f\sim g$ where $g$ is as above but with $c(x)$ replaced by $c$. Thus $g$ is differentiable with 
\[ g^{\prime}(x)=\frac{g(x)}{x}( \rho +\eps(x))\sim\rho \frac{g(x)}{x}. \]
In particular, every regularly varying function $f$ of index 1 is asymptotic to a differentiable function $g$ which satisfies $g^{\prime}(x)\sim\frac{g(x)}{x}$. 
See \cite{BGT} for a detailed exposition of the subject.\nl

\noindent
{\bf 3. Abel Functions}\nl
Given $f\in {\rm SIC}_\infty$ such that $f(x)>x$, the equation 
\[ F(f(x))=F(x)+1,\tag{A1}\]
called the {\em Abel Functional Equation}, always has infinitely many solutions $F\in {\rm SIC}_\infty$ (see \cite{KCG} or \cite{TH1}, section 3).  
As such, $F$ is called an {\em Abel function} of $f$. This is easily seen: define $F$ on $[A,f(A)]$ (with $A$ suitably large) such that $F$ is strictly increasing and continuous and $F(f(A))=F(A)+1$ and extend to $[A,\infty)$ via (A1). Furthermore, if $f\in {\rm D}^+_\infty$ then, by similar reasoning, we can produce (infinitely many) solutions $F\in {\rm D}^+_\infty$. Solutions of (A1) can be used to give fractional iterates of $f$, namely, via $f_\lambda :=  F^{-1}(F+\lambda)$. 

Now we show we can make solutions of (A1) satisfy $(R_0)$. \nl

\noindent
{\bf Proposition A.1}\nl
{\em Let $f\in {\rm SIC}_\infty$ be such that $f(x)\ge cx$ for some $c>1$ and $\frac{f(x)}{x}$ is increasing. Then every ${\rm SIC}_\infty$-Abel function $F$ of $f$ satisfies $(R_0)$. 

Furthermore, if $f\in {\rm D}^+_\infty$, then we may choose $F\in {\rm D}^+_\infty$, and as such $\frac{F(x)}{x}$ is decreasing.\nl

\noindent
Proof.}\, Since $\frac{f(x)}{x}$ is increasing, we have, for $\beta>\alpha$
\[ \frac{f^n(\beta)}{f^n(\alpha)} \ge \frac{f^{n-1}(\beta)}{f^{n-1}(\alpha)}\ge \cdots \ge \frac{\beta}{\alpha};\]
i.e. $\beta\le \alpha\frac{f^n(\beta)}{f^n(\alpha)}$. Let $y>0$ such that $x+y\le f(x)$. For $x$ sufficiently large, we may write $x=f^n(\alpha)$ and $x+y=f^n(\beta)$ for some $n\in\N$ and $\alpha\in [A,f(A)]$ and $\beta\in [A,f^2(A)]$. By above we note that $\beta\le\alpha + \frac{\alpha y}{x}$. Thus
\[ F(x+y)-F(x) = F(f^n(\beta))-F(f^n(\alpha)) = F(\beta)-F(\alpha)\le F\Bigl(\alpha + \frac{\alpha y}{x}\Bigr) - F(\alpha).\]
If now $y=o(x)$, then by continuity of $F$, the RHS above tends to 0. This shows $F$ satisfies $(R_0)$. 

Now suppose $f\in {\rm D}^+_\infty$, then we may choose $F\in {\rm D}^+_\infty$. As such, the above now gives
 \[ F(x+y)-F(x) \le C\frac{y}{x},\]
for some constant $C$. Hence
\[ \frac{F(x)}{x} - \frac{F(x+y)}{x+y} \ge F(x)\biggl(\frac{1}{x}-\frac{1}{x+y}\biggr) - \frac{Cy}{x(x+y)} = \frac{(F(x)-C)y}{x(x+y)}>0\]
for $x$ sufficiently large. Thus $\frac{F(x)}{x}$ is decreasing on each interval $[A,f(A)]$, which proves the result.  
\bo

Now we show that assumption (b) in the definition of a class generating sequence (prior to Definition 2.1) can, without loss of generality, be replaced by the stronger condition that $F_n^{\prime}$ itself is differentiable and $F_n^{\prime\prime}(x)\sim -F_n^{\prime}(x)/x$ (for $n\geq 2$). Further, in Proposition A.3 we show that such sequences can be generated easily. \nl

\noindent
{\bf Proposition A.2}\nl
{\em In the definition of a class generating sequence, $({\rm b})$ can be replaced by 

$({\rm b})^{\prime}$: $F_n=f_n^{-1}$ satisfies $(R_2)$ for $n=0,1$ and $(R_3)$ for $n\geq 2$, $F_n^{\prime}$ is differentiable and }
\[ F_n^{\prime\prime}(x)\sim -\frac{F_n^{\prime}(x)}{x}. \]

\noindent
{\em Proof.}\, By the above discussion, we have $F_n\sim 1/H_n$ where
\[ H_n(x) = cx\exp\Bigl\{ \int_A^x \frac{\eps(t)}{t}\, dt\Bigr\},   \]
where $\eps(x)\to 0$.

Now let $G_0=F_0$, $G_1=F_1$, and $G_n(x) = \int^x 1/H_n$ for $n\geq 2$. Then $G_n^{\prime}=1/H_n$ is regularly varying of index $-1$ (i.e. $G_n$ satisfies $(R_3)$ for $n\geq 2$.) Also $G_n^{\prime}\sim F_n^{\prime}$ and $G_n\sim F_n$. By Proposition 1.2(e), $\mathcal{B}_{F_n}=\mathcal{B}_{G_n}$ and $G_n$ satisfies $(R_0)$ for $n\geq 3$.  It remains to check $O_{G_{n+1}}(G_n)=-1$.

But for $n\geq 2$, $G_{n+1}(G_n)=G_{n+1}(F_n)+o(1)$. Since $G_{n+1}$ and $F_{n+1}$ are order-equivalent, $O_{G_{n+1}}(F_n)=  O_{F_{n+1}}(F_n)=-1$. Thus $O_{G_{n+1}}(G_n)=-1$. By Theorem 2.6, the classes generated by $\{F_n\}$ and $\{G_n\}$ are identical.
\bo

\noindent
{\bf Proposition A.3}\nl
{\em Let $f\in {\rm D}_{\infty}^+$ be such that $f^{\prime\prime}$ is continuous and satisfies 
\[ f^{\prime\prime}(x)\sim -\frac{f^{\prime}(x)}{x}.\]
(This implies $f$ satisfies $(R_3)$.) Then every $F\in {\rm D}_{\infty}^+$ solving $F=F(f)+1$ with $F^{\prime\prime}$ continuous satisfies
\[ F^{\prime\prime}(x)\sim -\frac{F^{\prime}(x)}{x} \]
(and hence satisfies $(R_3)$.) Furthermore, there exists at least one such $F$.}\nl

\noindent
{\em Proof.}\, First suppose there is an $F\in {\rm D}_{\infty}^+$ solving $F=F(f)+1$ with $F^{\prime\prime}$ continuous. Then $F(f)=F-1$, $F^{\prime}(f)f^{\prime}=F^{\prime}$ and
\[ F^{\prime\prime}(f)(f^{\prime})^2 + F^{\prime}(f)f^{\prime\prime} = F^{\prime\prime}. \]
Dividing through gives
\[ \frac{F^{\prime\prime}(f)f^{\prime}}{F^{\prime}(f)} + \frac{f^{\prime\prime}}{f^{\prime}} = \frac{F^{\prime\prime}}{F^{\prime}}. \]
Let $H(x) = -x\frac{F^{\prime\prime}(x)}{F^{\prime}(x)}-1$. We aim to show that $H(x)\to 0$ as $x\to\infty$. On rearrangeing the above we have
\[ H = \delta + \eta H(f), \tag{A2}\]
where $\delta$ and $\eta$ are the functions given by
\[ \eta(x) = \frac{xf^{\prime}(x)}{f(x)}, \qquad \delta(x)=1+\frac{xf^{\prime\prime}(x)}{f^{\prime}(x)} - \eta(x). \]
First we note that $\eta(x)\to 0$ as $x\to\infty$. For $f^{\prime\prime}<0$ (eventually), so $f^{\prime}$ decreases (to 0). Hence $f^{\prime}$ is regularly varying of index $-1$ and $f$ regularly varying of index 0 (see \cite{BGT}, p.59). Hence $\frac{xf^{\prime}(x)}{f(x)}\to 0$. By the assumption on $f$, $\delta(x)\to 0$ also.  

Let $\eps>0$ and choose $x_0$ large enough so that $|\delta(x)|<\eps$ and $|\eta(x)|<\alpha$ for $x\geq x_0$ (some $\alpha<1$). Then for $x\geq x_0$
\[ |H(x)| < \eps + \alpha|H(f(x))|. \]
By induction, 
\[ |H(f^{-n}(x)| <\eps(1+\alpha+\ldots +\alpha^{n-1}) + \alpha^n |H(x)|<\frac{\eps}{1-\alpha}+\alpha^n|H(x)|. \]
Every $y\geq x_0$ can be written uniquely as $y=f^{-n}(x)$ for some $n\in\N_0$ and $x\in [x_0,f^{-1}(x_0))$. Thus 
\[|H(y)|<\frac{\eps}{1-\alpha}+\alpha^n M, \]
where $M$ is the maximum value of $H$ on $[x_0,f^{-1}(x_0)]$. Thus $\limsup_{y\to\infty} |H(y)|\leq\frac{\eps}{1-\alpha}$. This holds for every $\eps>0$, so the limsup is zero; i.e. $H\to 0$, as required. 

Finally, to prove the existence of such an $F$, define $H$ continuous on some interval $[f(A),A]$ with $A$ suitably large such that (A2) holds at $A$. Extend $H$ to $[f(A),\infty)$ by (A2). Then define $F$ via
\[ \log F^\prime(x) =  - \int_A^x \frac{H(t)+1}{t}\, dt.\]
\bo

\noindent
{\bf 4. Hardy $L$-functions in classes}\nl
Let $\mathcal{H}$ denote the set of Hardy $L$-functions. We show that with $F_n=\Xi_n$, they  are contained in the classes $C_0$, $C_1$ or $C_2$. \nl

\noindent
{\bf Proposition A.4}\nl
{\em Let $f\in\mathcal{H}$ such that $f(x)\ge x+a$ for some $a>1$. Then $f\in C_0$, $C_1$ or $C_2$. Indeed, $f\in C_0^{(4)}$, $C_1^{(3)}$ or $C_2^{(2)}$.\nl

\noindent
Proof.}\, As we know, $O_\Xi(f)=k$ exists and is an integer. The assumption implies $k\ge 0$. If $k\ge 1$, then $O_F(f)=1$ where $F=\frac1k\Xi$. But $F^{-1}\in C_3^{(1)}$ since $O_{\Xi_4}(F^{-1})=1$. Thus $f\in C_2^{(2)}$ (since we already know that $f\in\mathcal{B}$). This leaves the case where $O_\Xi(f)=0$.

Now, as Hardy showed, there exists $n\in\Z$ and $\mu>0$ such that
\[  \log_n f(x)= (\log_n x)^{\mu+o(1)}.\tag{A3}\]
If $\mu>1$, then 
\[ \log_{n+2} f(x)= \log_{n+2} x + \log \mu+o(1).\]
(Note that $n+2\ge 0$ necessarily holds, otherwise $f(x)=x+o(1)$ follows.) Hence with $F(x)=\frac{\log_{n+2} x}{\log \mu}$, we have $O_F(f)=1$. If $n+2>0$ then $F^{-1}\in C_2^{(2)}$ (being a Hardy function of positive order). As such, $f\in C_1^{(3)}$. If $n=-2$, then $f(x)=x+\log\mu+o(1)$ is only possible if $\mu>e$. Then, as for the function $x+a$, $f\in C_0^{(3)}$. 

There remains the more tricky case $\mu=1$. Now $n\ge -1$ is forced. Replacing the $o(1)$ term in (A3) by $\frac{1}{h(x)}$ gives
\[ \log_{n+2} f(x)= \log_{n+2} x + \frac1{h(x)}.\]
Without loss of generality, we may assume $n$ is minimal; i.e. $1\prec h(x)\ll \log_{n+1} x$. (Here we use the fact that $h\in\mathcal{H}$.) We find $F$ such that $O_F(f)=1$ using Theorem 1.3. We have
\[  \Xi(f(x))-\Xi(x)=\Xi(\log_{n+2}f(x))-\Xi(\log_{n+2}x)\sim \frac1{h(x)\chi(\log_{n+2}x)}.\]
To apply Theorem 1.3, we need the LHS to be asymptotic to $\frac1{\chi(F(x))}$ for some $F\in\mathcal{B}$. This involves
\[  \chi(F(x))\sim h(x)\chi(\log_{n+2}x),\]
or, equivalently, using $\chi^{-1}(x)\sim \frac{x}{\chi(\log x)}$,
\[  F(x)\sim \frac{h(x)\chi(\log_{n+2}x)}{\chi(\log h(x) + \log\chi(\log_{n+2}x))}.\tag{A4}\]
We distinguish between the two possible cases: (i) $\log h(x)\succ \log_{n+3}x$ and (ii) $\log h(x)\sim c \log_{n+3}x$ for some\footnote{Here, if $c=0$, we mean $\log h(x)=o(\log_{n+3}x)$.} $c\ge 0$. If (ii) holds, then take $F(x)=\frac{h(x)}{c+1}\log_{n+2}x$, which is asymptotic to RHS (A4). As $F$ is a Hardy function, $F\in\mathcal{B}$, so by Theorem 1.3, $O_F(f)=1$. But $O_\Xi (F)=-(n+2)\le -1$, so from the first case, we know $F^{-1}\in C_2^{(2)}$. Hence $f\in C_1^{(3)}$. 

This leaves the case (i). As $\log_{n+3} x\prec \log h(x)\ll \log_{n+2}x$, we see that $O_\Xi(h)=-k$ where $k=n+1$ or $n+2$. So for large enough $r$, $\log_r h(x)\sim \log_{r+k}x$. The RHS of (A4) is therefore 
\[  \frac{h(x)\log_{n+2} x\cdots \log_{r+k-1}x\chi(\log_{r+k} x)}{\log h(x)\cdots \log_{r-1}h(x)\chi(\log_r h(x))}\sim
\frac{h(x)\log_{n+2} x\cdots \log_{r+k-1}x}{\log h(x)\cdots \log_{r-1}h(x)}.\]
The function on the right, being in $\mathcal{H}$, is in $\mathcal{B}$, so take this to be $F$. By Theorem 1.3, $O_F(f)=1$. But $O_\Xi(F)=O_\Xi(h)=-k$, so if $k>0$, we have $F^{-1}\in C_2^{(2)}$ and $f\in C_1^{(3)}$. The case $k=0$ can only occur when $n=-1$ and $f(x)\sim x$. But $f^n(\alpha)\ge an+\alpha$ for all $\alpha$ large enough, so
\[  F(f^n(\alpha))\sim n\le \frac1a f^n(\alpha).\]
Thus $F(x)\le bx$ for some $b<1$. For this case the above already shows that $F^{-1}$ lies in $C_1^{(3)}$ or $C_2^{(2)}$. Thus $f\in C_0^{(4)}$ or $C_1^{(3)}$.
\bo

We note that the same result holds for the class $E$ of Boshernitzan which, as mentioned in \S1, extends $\mathcal{H}$ to include solutions of certain algebraic differential equations and is closed under differentiation and integration. The growth properties of such functions are similar to those of $\mathcal{H}$; in particular (A3) holds. \nl

\noindent
{\bf 5. Construction of classes with less stringent regularity}\nl
In the article, we have developed a theory of classes which implicitly assumes a `best behaved' class generating sequence, by insisting all functions belong to $\mathcal{B}$. 
For example, if $F_1$ and $F_2$ are two ${\rm D}_{\infty}^+$-solutions of the Abel equation $F(\exp)=F+1$, then $F_1\sim F_2$. If $F_1,F_2\in \mathcal{B}$, then automatically $F_1^{\prime}\sim F_2^{\prime}$ which forces $F_1-F_2$ to be constant. Thus membership of $\mathcal{B}$ forces {\em uniqueness} (up to a constant) of this and any other Abel equation, and hence also forces unique fractional iterates. In general, for $F_1,F_2\not\in\mathcal{B}$, we only have $F_1^{\prime}\asymp F_2^{\prime}$.   
\nl

Below, we shall briefly describe the theory in the case when there is no (or no known) `best behaved' sequence. The equivalent notion of $\mathcal{B}_F$ is the following. 
\nl

\noindent
{\bf Definition}\nl
For $F\in {\rm D}_{\infty}^+$, let 
\[  \mathcal{B}_F^{\prime}=\{f\in {\rm D}_{\infty}^+ : (F\circ f)^{\prime}\asymp F^{\prime}\}.  \]
Writing $F^{\prime}=1/L$, this is equivalent to
\[ f^{\prime}\asymp \frac{L(f)}{L}. \]
Proposition 1.2 now becomes:\nl

\noindent
{\bf Proposition A.5}
\begin{enumerate}
\item $\mathcal{B}_F^{\prime}$ {\em is a group under composition}.
\item {\em If $F$ satisfies ($R_2$) and $f,g\in\mathcal{B}_F^{\prime}$, then $f\sim g$ implies $f^{\prime}\asymp g^{\prime}$. Under the stronger condition that $F$ satisfies $(R_3)$, then $f\asymp g$ implies $f^{\prime}\asymp g^{\prime}$. The special case $f(x)=o(x)$ implies $f^{\prime}\to 0$ also holds.} 
\item {\em Let $G\in {\rm D}_{\infty}^+$ satisfy $(R_3)$ and suppose $F\in\mathcal{B}_G^{\prime}$. Then
\begin{enumerate}
\item $\mathcal{B}_F^{\prime}\subset \mathcal{B}_G^{\prime}$; (i.e. $f\in\mathcal{B}_F^{\prime}\implies f\in\mathcal{B}_G^{\prime}$).
\item If $F(f)\asymp F$, then $f\in\mathcal{B}_G^{\prime}\implies f\in\mathcal{B}_F^{\prime}$.
\end{enumerate}
Thus $\mathcal{B}_F^{\prime}$ is the subgroup of functions $f$ in $\mathcal{B}_G^{\prime}$ for which $F\circ f\asymp F$.
\item If $F$ satisfies ($R_3$), then $f\in\mathcal{B}_F^{\prime}$ implies $\alpha f\in\mathcal{B}_F^{\prime}$ for every $\alpha>0$.
\item Let $F,G\in {\rm D}_{\infty}^+$. If $F^{\prime}\asymp G^{\prime}$, then $\mathcal{B}_F^{\prime}=\mathcal{B}_G^{\prime}$.}
\end{enumerate}

The proofs are analogous to those of Proposition 1.2, with $\asymp$ taking the place of $\sim$ in the appropriate places. We omit the details.\nl

The definition of the classes remains the same except that we now have $f_n\in\mathcal{B}_{F_{n+1}}^{\prime}$ in part (c). Letting
\[ \mathcal{B}^{\prime} = \bigcup_{n=0}^{\infty}\mathcal{B}_{F_n}^{\prime}, \]
Definition 2.1 becomes: for $n\geq 0$, $D_n^{(0)} = \{ f_n\}$, and supposing $D_n^{(k)}$ has been defined for all $n\geq 0$ and some $k\geq 0$, let 
\[ D_n^{(k+1)} = \{f\in\mathcal{B}^{\prime}:O_F(f) = 1\quad\mbox{for some $F$ with $F^{-1}\in D_{n+1}^{(k)}$}\}. \]
Now define classes $D_n$ by 
\[ D_n = \bigcup_{k=0}^{\infty} D_n^{(k)}. \]

Theorem 2.3, Corollary 2.4 and Theorem 2.6 generalise directly with $\mathcal{B}$ replaced by $\mathcal{B}^{\prime}$.
By suitably adapting the choice of $h$ in Theorem 2.7, this theorem and consequently Theorem 2.8 also hold in this more general setting. Basically, one shows that there exists functions in $C_n^{(k+1)}$ which are larger (resp. smaller) than all functions in $D_n^{(k)}$. For this we require the following generalization of Theorem 1.3:\nl

\noindent
{\bf Theorem A.6}\nl 
{\em Let $F\in\mathcal{B}_G^{\prime}$ for some function $G$ satisfying $(R_1)$. Then}
\[ F(f)-F\asymp 1 \Longleftrightarrow G(f) - G\asymp G^{\prime}(F). \]
{\em Proof.}\, Write $G(F)=h(G)$ where $h^{\prime}\asymp 1$. 

$(\Leftarrow)$\, We have
\[ G(F(f)) = h(G(f)) = h(G+\theta G^{\prime}(F)) = h(G)+\theta_1 G^{\prime}(F) = G(F)+\theta_1 G^{\prime}(F) = G(F+\theta_2) \]
where $\theta,\theta_1,\theta_2\asymp 1$. Thus $F(f)-F\asymp 1$.\nl

$(\Rightarrow)$\, We have
\[ G(F) = h^{-1}(G(F(f))) = h^{-1}(G(F+\theta)) = h^{-1}(G(F)+\theta_1 G^{\prime}(F)) = h^{-1}(G(F))+\theta_2 G^{\prime}(F) = G+\theta_2 G^{\prime}(F), \]
where again, $\theta,\theta_1,\theta_2$ are all $\asymp 1$ (though different from before). Thus $G(f)-G\asymp G^{\prime}(F)$.
\bo

\end{document}